\input amstex
\documentstyle{amsppt}

\document
\vsize=6 in
\hsize=5.5 in
\magnification=1200

\topmatter
\title  A criterion for positive polynomials \endtitle
\author Fernando Cukierman \endauthor
\date{June 2004}\enddate
\thanks 
{I warmly thank Emili Bifet, Alicia Dickenstein, Gabriel Taubin and Bernd Sturmfels 
for several communications. This work was presented at the S\'eminaire de 
G\'eom\'etrie Analitique of the Universit\'e de Rennes;
I also thank Jorge Vitorio Pereira,  Dominique Cerveau, Frank Loray and Michel Coste
for their hospitality, interest and useful comments.}
\endthanks
\address {Departamento de Matem\'atica, Ciudad Universitaria,
(1428) Buenos Aires, Argentina } \endaddress
\email {fcukier\@dm.uba.ar} \endemail
 \endkeywords
  \endsubjclass
\abstract {Let $F$ be a homogeneous real polynomial of even degree in any 
number of variables.  We consider the problem of giving explicit conditions
on the coefficients so that $F$  is positive definite or positive semi-definite. 
In this note we produce a necessary condition for positivity and a sufficient
condition for non-negativity, in terms of positivity or semi-positivity
of a one-variable {\it characteristic polynomial} of $F$.  Also, we review
another well-known sufficient condition. } \endabstract
\endtopmatter

\comment
\NoRunningHeads
\NoPageNumbers
\leftheadtext{  }
\rightheadtext{  }
\endcomment

\define \R{\Bbb R}
\define \C{\Bbb C}
\define \N{\Bbb N}

\flushpar
{\bf \S 1}
\newline
\newline
(1.1) Let $F$ be a homogeneous polynomial of degree $d$ in $n$ 
variables $x_1, \dots, x_n$
with coefficients in a field $K$.
We will denote 
$$K(n, d)$$
the $K$-vector space of all such polynomials. Its dimension is $N = \binom {n-1+d}{d}$.
\newline
\newline
For $K = \R$, the field of real numbers,
we shall say that $F$ is {\bf positive} (resp. {\bf non-negative}), 
written $F > 0$ (resp. $F \ge 0$), 
if $F(x) > 0$ for all $x \in \R^n -\{0\}$ (resp. $F(x) \ge 0$ for all $x \in \R^n$).
\newline
\newline
We are interested in obtaining conditions on the coefficients of $F$ equivalent to 
$F > 0$ or to $F \ge 0$. We assume $d$ is even, so that positive polynomials
exist.
\newline
\newline
(1.2) In case $d=2$ such conditions are given by the well-known Sylvester's criterion:
\newline
If $F(x) = \sum_{1 \le i, j \le n}  F_{ij}  x_i  x_j$  (with $F_{ij} = F_{ji} \in \R$) then
$$F > 0  \  \text{ if and only if }  \  D_r(F) > 0,  \  r = 1, \dots, n$$
where $D_r(F) = \text{det } (F_{ij})_{1 \le i, j \le r}$ is the
$r \times r$ principal minor of the $n \times n$ symmetric matrix representing $F$.
\newline
\newline
Let us remark that the conditions $F \ge 0$ and $D_r(F) \ge 0,  \  r = 1, \dots, n$
are not equivalent. Actually,  $F \ge 0$ is equivalent to 
$$D_J(F) \ge 0   \text { for all }  J \subset \{1, \dots, n\}$$
where $D_J(F) = \text{det } (F_{ij})_{i, j \in J}$ (see \cite {\bf  Ga}).
\newline
\newline
(1.3) For the case $n=2$ of binary forms we can use Sylvester's formulation
of Sturm's Theorem (\cite {\bf  Ga}, \cite {\bf  De}, \cite {\bf  P}). To recall this,
let $p \in \R[t]$ be a (monic) polynomial of 
degree $d$ in one variable, over the real numbers. Consider the finite dimensional
$\R$-algebra $A = \R[t] / (p)$ and denote its trace linear form  $\text{tr}: A \to \R$. 
For each $u \in A$ define a quadratic form on $A$ by 
$$Q_u(x) = \text{tr}(u x^2)$$
Let $R(p) \subset \C$ be the set of complex roots of $p$. 
Denoting sg the signature of a quadratic form, Sylvester's theorem asserts that
$$
\align
\text{rank } (Q_u) &= \text{card } \{r \in R(p) /   \ u(r) \ne 0 \} \\ 
\text{sg}(Q_u) &= \text{card } \{r \in R(p) \cap \R / \  u(r)  >  0 \} - 
\text{card } \{r \in R(p) \cap \R /  \ u(r) <  0 \}
\endalign
$$ 
Suppose $p(0) \ne 0$ and denote $P$ (resp. $N$) the number of positive 
(resp. negative) real roots of $p$.
Choosing $u = 1$, we have that $\text{sg}(Q_1)$ is the
number $P + N$ of real roots of $p$. 
For $u = t$ we obtain $ \text{sg}(Q_t)  = P - N$.
Hence, $2 P = \text{sg}(Q_1) + \text{sg}(Q_t)$. In particular,
$P = 0$ (i. e. $p(t) > 0$ for $t > 0$) if and only if 
$\text{sg}(Q_1) + \text{sg}(Q_t) = 0$, a fact that will be useful later.
\comment
como se calcula la signatura de una forma cuadratica en general ?
creo es numero de cambios de signos de dets principales
\endcomment
\newline
\newline
(1.4) Let us denote $\nabla(n, d, \C) \subset \C(n, d)$ the set of 
singular polynomials of degree $d$ in $n$ variables, over the complex numbers. 
That is,
$$\nabla(n, d, \C) = \{ F \in  \C(n, d) \  / \ \frac {\partial F}{\partial x_i} (x) = 0,
\forall i, \text { for some } x \in \C^n-\{0\} \}$$
It is known (see e.g. \cite {\bf  GKZ}) that $\nabla(n, d, \C)$ is an irreducible algebraic 
hypersurface of degree 
$$D = n (d-1)^{n-1}$$ 
defined over the rational numbers. Therefore, there exists a polynomial 
(unique up to multiplicative constant)
$$\Delta = \Delta(n, d)$$
called the {\bf discriminant}, such that 
$$\nabla(n, d, \C) = \{F \in  \C(n, d) \ / \  \Delta(F) = 0\}$$ 
More precisely, writing a general polynomial 
$$F = \sum_{|\lambda| = d} F_{\lambda} \  x^{\lambda} \ \   \in \ \  K(n, d)$$
where $\lambda = (\lambda_1, \dots, \lambda_n) \in \Bbb N^n$, 
$|\lambda| = \sum_i \lambda_i$, $x^{\lambda} = \prod_i x_i^{\lambda_i}$ 
and $F_{\lambda} \in K$, so that  $T_{\lambda}(F) = F_{\lambda}$ 
are coordinate functions
on the vector space $K(n, d)$, we know $\Delta$ is a homogeneous polynomial 
in variables  $T_{\lambda}$, of degree $D = n (d-1)^{n-1}$ 
and with rational coefficients.
In other terms, $\Delta$ is an element of the $D$-th symmetric power of
the rational vector space dual to $\Bbb Q(n, d)$.
\newline
\newline
We shall normalize $\Delta$ so that $\Delta(J) = 1$ where 
$J \in \C(n, d) - \nabla(n, d, \C)$ is the polynomial
$$J(x) = \sum_{1 \le j \le n} x_j^d$$
(1.5) Restricting to the real numbers, we denote
$$
\align
\nabla = \nabla(n, d, \R) &= \nabla(n, d, \C) \cap \R(n, d) \\
&= \{ F \in  \R(n, d) \  /  \ \frac {\partial F}{\partial x_i} (x) = 0, \forall i,
\ \text {for some } x \in \C^n, x \ne 0 \} \\
&= \{F \in \R(n, d) \  / \  \Delta(F) = 0\}
\endalign
$$
the set of real polynomials which have a singular point, real or complex.
\newline
\newline
(1.6) Let us denote
$$\Cal P = \Cal P(n, d) = \{ F \in  \R(n, d) \  / \  F > 0 \}$$
the set of all positive polynomials. It is easy to verify that $\Cal P$ 
is an open convex cone in the vector space $\R(n, d)$.
Here by "cone" we mean a set which is stable under multiplication by $\R_{>0}$.
It is also easy to see that the closure of $\Cal P$ (in the usual topology of $\R(n, d)$)
 is the closed convex cone
$$\bar\Cal P = {\bar \Cal P}(n, d) = \{ F \in  \R(n, d) \ /  \ F \ge 0 \}$$
With slight abuse of notation we write $\Delta: \R(n, d)  \to \R$ for the polynomial 
function induced by the polynomial $\Delta$.
\newline
\newline
(1.7) {\bf Theorem:} $\Cal P - \nabla$ is connected.
\newline
\newline
{\bf Proof:} Let $F \in \Cal P \cap \nabla$ be a positive singular polynomial.
Let $x \in \C^n$ be a singular point of $F$. Since $F(x) = 0$ and $F > 0$ it
follows that $x \in \C^n -\R^n$. The complex conjugate $\bar  x \in \C^n$
is also a singular point of $F$, because $F$ has real coefficients.
It follows that $F$ has two distinct singular points in $\C^n$.
\newline
\newline
The idea now is that polynomials with at least two singular points occur 
in codimension at least two and hence they do not disconnect $\Cal P$.
To prove this, let us denote
$$\nabla_2(n, d, \C) =  \text{closure} \{ F \in  \C(n, d) \  / \ 
\frac {\partial F}{\partial x_i} (x) = \frac {\partial F}{\partial x_i} (y) = 0, \forall i, 
\text { for some } x \ne y \in \C^n - \{0\}  \}$$
the Zariski closure of the set of complex polynomials with at least two singular points.
By the usual incidence correspondence argument,
$\nabla_2(n, d, \C) \subset  \C(n, d)$ is a complex algebraic variety
of codimension two. Let us denote its real points
$$\nabla_2 = \nabla_2(n, d, \R) = \nabla_2(n, d, \C) \cap \R(n, d)$$
We claim that $\Cal P \cap \nabla = \Cal P \cap \nabla_2$. The inclusion $\subset$
was observed just above, while the other one is clear since $\nabla_2 \subset \nabla$.
\newline
\newline
Then $\Cal P \cap \nabla   \subset \R(n, d)$ is an open subset of
a real algebraic variety $\nabla_2$ of real dimension $N - 2$, where
$N = \binom {n-1+d}{d} = \text {dim } \R(n, d)$.  
\newline
\newline
The Theorem now follows from (1.8) below,
which is surely a well-known statement, but I shall give a proof
because I could not find a suitable reference.
\newline
\newline
(1.8) {\bf Proposition:} Let $P  \subset \R^N$ be a connected open set,
$Y \subset \R^N$ a real algebraic variety of dimension $d$ and denote
$X = P \cap Y$. Then
\newline
\newline
a) For any family of supports $h$ and every sheaf $\Cal L$
of abelian groups, $H^j_{h}(X, \Cal L) = 0$ for all $j > d$
(i. e. the $h$-cohomological dimension of $X$ is $\le d$). 
\newline
\newline
In particular, $H_c^{j}(X, \Bbb Z) = 0$ for $j > d$,
where $H_c$ denotes cohomology with compact supports.
\newline
\newline
b) If  $d \le N-2$ then $P - X$ is connected.
\newline
\newline
{\bf Proof:} We refer to \cite {\bf  Go} for general definitions.
\newline
To prove a), let us denote $S \subset Y$ the set of singular points,
$T \subset Y$ the union of the irreducible components of $Y$
of dimension $< d$ and $A = P \cap (S \cup T)$, which is
a closed subset of $X$. Let us remark that 
\newline
i) $S \cup T \subset \R^N$ is a real algebraic variety of dimension $< d$.
\newline
ii) $Y - (S \cup T)$ is a smooth manifold of dimension $d$.
Then, its open subset $X - A$ is also a smooth manifold of dimension $d$.
\newline
Now we apply to $(X, A)$ the theory of \cite {\bf  Go}, (4.10). From the exact 
sequence of sheaves on $X$
$$0 \to \Cal L_{X-A} \to \Cal L \to \Cal L_{A} \to 0$$
one obtains a long exact sequence \cite {\bf  Go}, (4.10.1) of cohomology with 
supports in $h$
$$ \dots \to H^j_{h | X-A}(X - A, \Cal L) \to H^j_{h}(X, \Cal L) 
\to H^j_{h | A}(A, \Cal L) \to H^{j+1}_{h | X -A}(X - A, \Cal L) \to \dots$$
Let $j > d$. By induction on $d$ and due to i) we know that  
$H^j_{h | A}(A, \Cal L) = 0$.
Since a smooth manifold of dimension $d$ has $h$-cohomological dimension $d$
(\cite {\bf  Go}, (4.14.1) and (5.12)) we have $ H^j_{h | X - A}(X - A, \Cal L) = 0$.
From the exact sequence it follows that $H^j_{h}(X, \Cal L) = 0$, as claimed.
\newline
In particular, when $\Cal L = \Bbb Z_X$ is the constant
sheaf $\Bbb Z$ on $X$ and $h$ is the family of compact subsets of $X$,
we obtain the statement about cohomology with compact supports.
\newline
\newline
To prove b) we apply the Poincare-Lefschetz
duality theorem (\cite {\bf  Do}, (7.13)) for the $N$-dimensional manifold $P$
and the closed set $X \subset P$. This Theorem implies that
$$H_c^{N-1}(X) = H_1(P, P - X)$$
Combining with a) we obtain $H_1(P, P - X) = 0$ and from the homology sequence
of the pair $(P, P - X)$ we deduce
$H_0(P - X) = H_0(P) = \Bbb Z$, Q.E.D.
\newline
\newline
Now we deduce a Corollary of importance for our present purpose.
\newline
\newline
(1.9) {\bf Corollary:} Let $F \in \R(n, d)$. If $F \ge 0$ then $\Delta(F) \ge 0$.
\newline
\newline
{\bf Proof:} We want to show that $\Delta$ is non-negative on $\bar\Cal P$.
It is clear that 
$$\Cal P - \nabla \subset (\Delta > 0) \cup (\Delta < 0)$$
where we denote $(\Delta > 0) = \{F / \  \Delta(F) > 0 \}$.
By (1.7), $\Cal P - \nabla$ is connected and the polynomial $J$
defined in (1.4) above belongs to the intersection $(\Cal P - \nabla) \cap (\Delta > 0)$,
so we obtain $\Cal P - \nabla \subset (\Delta > 0)$, that is,
$\Cal P  \subset (\Delta \ge 0)$, and by continuity
$\bar\Cal P  \subset (\Delta \ge 0)$, as we wanted to prove.
\newline
\newline
(1.10) {\bf Example:} Les us consider the case $(n, d) = (2, 4)$
of binary quartics, written 
$F = \sum_{0 \le i \le 4} F_i  \ x_1^i  x_2^{4-i}$.
Over the complex numbers we have
 $\nabla(2, 4, \C) \subset \C(2, 4)$ 
of respective complex dimensions 4 and 5. Their real parts
$\nabla \subset \R(2, 4)$ have real dimensions 4 and 5, but
$\nabla$ contains the open set $\Cal P \cap \nabla$ which has 
dimension 3
and is a "component" of $\nabla$ in the sense that it is not
contained in the closure of $\nabla - \Cal P  \cap \nabla$,
with respect to the usual (not Zariski) topology of $\R(2, 4)$.

More explicitly, $\Cal P \cap \nabla$
consists of the quartics with two double complex roots, i. e. of the form
$F = (a x_1^2 + b x_1 x_2 + c x_2^2)^2$ with $b^2 - 4 a c  <  0$
(so we see again it has dimension 3) and it is clear that these 
are not limit of real quartics with only one double root.

Let $V \subset \R(2, 4)$ denote the 3-dimensional subspace of those $F$'s 
with $F_4 =1$ and  $F_3=0$. 
We refer to \cite {\bf  GKZ}, page 381, for a drawing of the two-dimensional
real algebraic variety $\nabla \cap V$, but let's point out that
$(\Cal P \cap \nabla) \cap V$
(curve of quartics of the form $(x_1^2 + c x_2^2)^2$ with $c > 0$)
is lacking in that picture and should be added as a curve pointing out
and going through the point $F=x_1^4$ marked "quadruple root".
\newline
\newline
(1.11) {\bf Definition:} For  $F \in K(n, d)$
we define the {\bf characteristic polynomial} of $F$ (with respect to $J$)
$$\chi(F; J)(t) = \Delta(F + t J) \  \in \ K[t]$$
where $J(x) = \sum_{1 \le j \le n} x_j^d \ $, as in (1.4). 
\newline
\newline
(1.12) {\bf Remarks:} 
\newline
a) The definition depends on the choice of $J$, but for simplicity we may
write $\chi(F)$ instead of $\chi(F; J)$. In fact, our choice of a positive $J$ is rather arbitrary.
\newline
b) Since $J$ and $\Delta$ have rational coefficients it follows that $\chi(F)$ has 
coefficients in $K$ if $F$ has coefficients in $K$.
\newline
c) $\chi(F)$ is a polynomial in $t$ of degree $D$ as in (1.4), and we may write
$$\chi(F)(t) = \sum_{0 \le j \le D} \Delta_j(F) \  t^j $$
where $\Delta_0 = \Delta$ and $\Delta_j$ is a homogeneous polynomial 
in the coefficients of $F$, of degree $D - j$ for $j = 0 , \dots, D$.
Also, by our normalization of  $\Delta$ 
(see (1.4)) it follows that $\chi(F)$ is monic.
\newline
d) The roots of $\chi(F)$ are the values of $t$ such that $F + t J$
is singular, that is,  they parametrize the intersections of the discriminant hypersurface 
$\nabla$ with the pencil spanned by $F$ and $J$.
These roots may have the right to be called "eigenvalues of $F$" (with respect to $J$).
\newline
\newline
The next Proposition gives a necessary condition for non-negativity.
\newline
\newline
(1.13) {\bf Proposition:} Let $F \in \R(n, d)$. If $F \ge 0$ then 
$\chi(F)(t) \ge 0$ for all $t \in \R_{\ge 0}$.
\newline
\newline
{\bf Proof:} Since $J \ge 0$, for $F \ge 0$  and $t \ge 0$ we have $F + t J \ge 0$.
By Corollary (1.9), $\chi(F)(t) = \Delta(F + t J) \ge 0$, as wanted.
\newline
\newline
(1.14) {\bf Remarks:} 
\newline
a) The Proposition gives a necessary condition for non-negativity, 
in principle explicitly computable  via applying the Sylvester criterion (1.3)
to the characteristic polynomial $\chi(F)$.
\comment
ver comentario anterior sobre >0 y \ge 0
\endcomment
\newline
b) As kindly pointed out to us by Jiawang Nie (Berkeley), 
the converse to (1.13) fails for $d=2$.
\newline
c) Regarding the choice of $J$ mentioned in (1.12)a),
let us fix a finite number of positive polynomials $J_1, \dots, J_m \in \R(n, d)$.
As in (1.13), if $F \ge 0$ then 
$\chi(F; J_i)(t) \ge 0$ for all $t \in \R_{\ge 0}$ and all $i = 1, \dots, m$.
It would be interesting to know if this necessary condition for non-negativity
is also sufficient, for some choice of $m$ and $J_1, \dots, J_m$.
\newline
d) Similarly, for $F \in \R(n, d)$ define a generalized characteristic polynomial
$$\chi(F; J_1, \dots, J_m)(t_1, \dots, t_m) = \Delta(F + \sum_{1=1}^m t_i J_i) 
\  \in \R[t_1, \dots, t_m]$$
The same proof of (1.13) gives that a necessary condition for $F \ge 0$ is
$$\chi(F; J_1, \dots, J_m)(t_1, \dots, t_m) \ge 0,  \ \text{ for } t_i \ge 0$$
\newline
As kindly pointed out to us by Michel Coste, the converse to (1.13) holds
replacing $\ge$ by $>$. More precisely,
\newline
\newline
(1.15) {\bf Proposition:} Let $F \in \R(n, d)$. If $\chi(F)(t) > 0$ for all $t \in \R_{\ge 0}$ then $F > 0$.
\newline
\newline
{\bf Proof:}  If $\chi(F)(t) = \Delta(F + t J) > 0$ for all $t \in \R_{\ge 0}$ then the family of
polynomials $F + t J$ with $t \ge 0$ has constant topological type.
Then, since the set of zeros of $J$ is equal to $\{0\}$, the same
is true for $F + t J$ for all $t \ge 0$. In particular, the set of zeros of $F$ is $\{0\}$, so that
$F > 0$ or $F < 0$. If  $F < 0$ then there exists  $t \ge 0$ such that $F + t J$ has a non-trivial zero,
which is a contradiction. Therefore $F > 0$, as wanted.
\comment
\newline
\newline
(1.16) {\bf Remark:} (Michel Coste) Denote 
$X =\Cal P \cap \nabla = \Cal P \cap \nabla_2  \subset  \R(n, d)$
as in the Proof of (1.7). For $J \in \R(n, d)$ denote $H_J \subset  \R(n, d)$
the cone over $X$ with vertex $J$. Since $X$ has codimension two 
(see Proof of (1.7)) $H_J$ is a hypersurface.  
Choose $m$ and $J_1, \dots, J_m \in  \Cal P $ such that
$ H_{J_1}  \cap \dots \cap H_{J_m} = X$. Then, if
$\chi(F; J_i)(t) \ge 0,  \ \text{ for } t \ge 0$ for all  
$i = 1, \dots, m$ it follows that $F \ge 0$.
This raises the problem of finding such $m$ and $J_1, \dots, J_m$.
Is is true that I can take $m=2$ ?  For any projective variety ?
\endcomment
\newline
\newline
Now we obtain further conditions for positivity by
combining the previous constructions with the operation
of restriction to a linear subspace.
\newline
\newline
(1.16) Let $V$ be a finite dimensional real vector space 
and denote ${S^d}(V^*)$ the $d$-th symmetric power of the dual of $V$,
thought of as the space of homogeneous polynomials of degree $d$ in $V$.
As in (1.4) denote $\nabla_V \subset S^d(V^*)$ the set of 
singular polynomials and $\Delta_V \in S^D(S^d(V))$ the discriminant.
\newline
\newline
Let $V \subset \R^n$ be a linear subspace.
For $F \in \R(n, d) = S^d(({\R^n})^*)$, denote $F_V \in S^d(V^*)$ 
the restriction of $F$ to $V$. We may now formulate a stronger
necessary condition for non-negativity:
\newline
\newline
(1.17) {\bf Proposition:} Let $F \in \R(n, d)$.
If $F \ge 0$ in $\R^n$  then  $\Delta_V (F_V) \ge 0$  for 
every linear subspace  $V \subset \R^n $. Also, the characteristic
polynomial $\chi_V(F)(t) = \Delta_V(F_V + t J_V)$ is  $\ge 0$ for $t \ge 0$
\newline
\newline
{\bf Proof:} It is clear that if $F \ge 0$ in $\R^n$ then its restriction
$F_V$ is also $\ge 0$ in $V$. Applying (1.9) in $V$ we obtain 
$\Delta_V (F_V) \ge 0$, as claimed. The claim about $\chi_V(F)$
in immediate as in (1.13).
\newline
\newline
(1.18) {\bf Example:} For  $J \subset \{1, \dots, n\}$
let $V_J = \{ x \in \R^n / x_i = 0 \text{ for } i \notin J \}$,
the $J$-th coordinate plane. 
It follows that if $F \ge 0$ in $\R^n$ then 
$$\Delta_{V_J}(F) \ge 0,  \text{ for all } J \subset \{1, \dots, n\}$$
Notice that in case $d=2$, these $\Delta_{V_J}$ coincide with the diagonal 
minors $D_J$ of (1.2). 
\newline
\newline
(1.19) {\bf Remark:} In case $d=2$, the conditions of (1.18) 
are equivalent to $F \ge 0$ (see (1.2)). For $d > 2$ this is no longer true,
as it is easily seen in the case of binary forms.
\newline
\comment
quiza sea cierto algo del tipo siguiente:
existe un numero finito de subespacios V_i tales que F \ge 0 es
equivalente a las condiciones
$\Delta_{V_i} (F_V) \ge 0$ para todo i
$F_V_j \ge 0$ para los j tales que dim V_j = 2 (formas binarias, 
para las cuales se usa el criterio de sylvester
\endcomment
\newline
\newline
\flushpar
{\bf \S 2}
\newline
\newline
In this section we review a well-known sufficient condition for non-negativity of 
$F \in R(n, d)$ in terms of a quadratic form $h(F)$ associated to 
$F$ (see \cite {\bf R}). 
\newline
\newline
(2.1) Let $V$ be a vector space of dimension $n$ over a field $K$.
The multiplication of the symmetric algebra $S(V^*)$ induces maps
$$\mu: S^m(S^d(V^*)) \to S^{md}(V^*)$$
These are homomorphisms between linear representations of $GL(V)$.
They may be interpreted geometrically as the pull-back of
homogeneous polynomials of degree $m$ under the $d$-th Veronese map
$$\Bbb P V \to \Bbb P S^d(V)$$
sending $v \in V $ to $ v^d \in S^d(V)$  (see \cite {\bf FH}).
\newline
\newline
Let $K=\R$ and fix $m, d$. Consider the vector space $U = S^d(V)$
and suppose $G \in S^m(U^*)$ is non-negative (resp. positive).
Then the restriction of $G$ to the $d$-th Veronese variety in $\Bbb P U$
is clearly non-negative (resp. positive) and hence $\mu(G) \in S^{md}(V^*)$ 
is non-negative (resp. positive). 
\newline
\newline
For $m=2$ in particular, we have 
$$\mu: S^2(S^d(V^*)) \to S^{2d}(V^*)$$ 
and if the quadratic form $G \in S^2(U^*)$ is $\ge 0$ (resp. $>0$)
then the homogeneous polynomial $\mu(G) \in S^{2d}(V^*)$ is 
$\ge 0$ (resp. $>0$).
\newline
\newline
(2.2) On the other hand, suppose we have a map 
$$h: S^{2d}(V^*) \to S^2(S^d(V^*))$$ 
such that $\mu \circ h = \text{identity}$. It follows  that  
if the quadratic form $h(F)$ on $U = S^d(V)$ is $\ge 0$ (resp. $>0$) then
$F = \mu(h(F)) \in S^{2d}(V^*)$ is $\ge 0$ (resp. $>0$). 
This is the sufficient condition mentioned above. 
\newline
\newline
(2.3) What we shall do next is to explicitly construct such an $h$.
It will be the well-known Hankel quadratic form $h(F)$ associated 
to a homogeneous polynomial $F \in S^{2d}(V^*)$ of even degree $2d$ 
(see \cite {\bf R}). Our construction will be based on the co-algebra structure
of the symmetric algebra. 
It will follow in particular that $h$ is linear and $GL(V)$-equivariant. 
The equivariance does not
seem obvious from the construction in  \cite {\bf R}, which is based on an inner 
product on $S^{2d}(V^*)$. Let us also remark that since $S^{2d}(V^*)$
is an irreducible representation of $GL(V)$, such an equivariant $h$ is unique
up to multiplicative constant. Thus, the construction below may be considered
as another example of plethysm as in \cite {\bf FH}. 
\comment
For its interpretation as
geometric plethysm see ( , ).  HACERLO. 
creo tiene que ver con la polaridad respecto a la curva de veronese,
ver telling, en el caso d=4
\endcomment
\newline
\newline
(2.4) To start with the construction, let $V$ be a vector space of dimension 
$n$ over a field $K$. For each $d \in \N$ we have a natural map
$$V^{\otimes d} \otimes {V^*}^{\otimes d} \to K$$
given on elementary tensors by
$$(v_1 \otimes \dots \otimes v_d) \otimes 
(\varphi_1 \otimes \dots \otimes \varphi_d) \mapsto
\frac{1}{d !} \sum_{\sigma \in {\Bbb S}_d} \prod_{i=1}^d <\varphi_i, v_{\sigma(i)}>$$
This map factors through the quotient and gives a map of symmetric powers
$$( , ): S^d(V) \otimes S^d(V^*) \to K$$
with similar formula for elementary tensors (monomials).
This induces a linear map (called "polarization", see e. g. \cite {\bf D})
$$\wp: S^d(V) \to (S^d(V^*))^*$$
Let us remark that $(, ) $ and $\wp$ are equivariant for the natural
actions of $GL(V)$.
\newline
\newline
Let $\{e_1, \dots, e_n\}$ be an ordered basis of $V$ and denote
$\{x_1, \dots, x_n\}$ the dual basis of $V^*$, so that $<x_i, e_j> = \delta_{ij}$.
For $\alpha = (\alpha_1, \dots, \alpha_n) \in \N^n$ with $|\alpha| = \sum_i \alpha_i = d$,
denote as usual 
$e^{\alpha} = e_1^{\alpha_1} e_2^{\alpha_2} \dots e_n^{\alpha_n} \in S^d(V)$ and
$x^{\alpha} = x_1^{\alpha_1} x_2^{\alpha_2} \dots x_n^{\alpha_n} \in S^d(V^*)$.
Then $\{e^{\alpha}\}_{|\alpha| = d}$ (resp. $\{x^{\alpha}\}_{|\alpha| = d}$)
is a basis of $S^d(V)$ (resp. of  ${S^d}(V^*)$). Also, it is easy to check
from the explicit formulas above that
$$(x^{\beta}, e^{\alpha}) = \frac{\alpha !}{d!} \ \delta_{\alpha \beta}$$
where $\alpha ! = \prod_i \alpha_i !$. It follows that
$$\wp(e^{\alpha})(\frac{d!}{\beta!} \ x^{\beta}) =  \delta_{\alpha \beta}$$
Writing $D^{\alpha} = \wp(e^{\alpha}) \in (S^d(V^*))^*$ we have
$$D^{\alpha}(\frac{d!}{\beta!} \ x^{\beta}) =  \delta_{\alpha \beta}$$ 
and the isomorphism $\wp$ may be written as
$$\wp(\sum_{|\alpha|=d} a_{\alpha} e^{\alpha}) = 
\sum_{|\alpha|=d} a_{\alpha} D^{\alpha}$$
(2.5) Now we look at the structure of co-algebra in the symmetric algebra
$$S(V) = \bigoplus_{n \in \N} S^n(V)$$ 
We consider the multiplication map
of the symmetric algebra $S(V^*)$
$$\mu: S^d(V^*) \otimes S^e(V^*) \to S^{d+e}(V^*)$$
and the diagram 
$$
\CD
(S^{d+e}(V^*))^* @>{\mu^*}>> (S^d(V^*) \otimes S^e(V^*))^* 
\cong (S^d(V^*))^* \otimes (S^e(V^*))^* \\
@A{\wp}AA   @AA{\wp \otimes \wp}A \\
S^{d+e}(V) @>{h}>> S^d(V) \otimes S^e(V) 
\endCD
$$
where the vertical arrows are isomorphisms and $h$ is defined so that the diagram commutes.
It easily follows from the definitions that the effect of $h$ on basis elements is
$$h(e^{\gamma}) = \sum \Sb |\alpha| = d \\  |\beta| = e \\ \alpha + \beta = \gamma \endSb
c_{\alpha \beta } \  e^{\alpha} \otimes e^{\beta}$$
where $c_{\alpha \beta} = \frac{d!}{\alpha!} \frac{e!}{\beta!} \frac{(\alpha+\beta)!}{(d+e)!}$.
In terms of the basis elements $E^{\alpha} = \frac{d!}{\alpha!} \ e^{\alpha} \in S^d(V)$,
$$h(E^{\gamma}) = \sum \Sb |\alpha| = d \\  |\beta| = e \\ \alpha + \beta = \gamma \endSb
 \  E^{\alpha} \otimes E^{\beta}$$
Applying this to $V^*$ we obtain $GL(V)$-equivariant maps
$$h: S^{d+e}(V^*) \to S^d(V^*) \otimes S^e(V^*)$$ 
such that for $X^{\alpha} = \frac{d!}{\alpha!} \ x^{\alpha} \in S^d(V^*)$
$$h(X^{\gamma}) = \sum \Sb |\alpha| = d \\  |\beta| = e \\ \alpha + \beta = \gamma \endSb
 \  X^{\alpha} \otimes X^{\beta}$$
(2.6) {\bf Proposition:} With the notation above and any $d, e$, the composition
$$S^{d+e}(V^*) @>{h}>> S^d(V^*) \otimes S^e(V^*) @>{\mu}>> S^{d+e}(V^*)$$
is the identity.
\newline
\newline
{\bf Proof:} Let us compute on basis elements $x^{\gamma}$ as above
$$
\mu(h(x^{\gamma})) = \mu(\sum \Sb |\alpha| = d \\  
|\beta| = e \\ \alpha + \beta = \gamma \endSb
c_{\alpha \beta } \  x^{\alpha} \otimes x^{\beta}) =
\sum \Sb |\alpha| = d \\  
|\beta| = e \\ \alpha + \beta = \gamma \endSb
c_{\alpha \beta } \  x^{\alpha + \beta} = 
(\sum \Sb |\alpha| = d \\  
|\beta| = e \\ \alpha + \beta = \gamma \endSb
c_{\alpha \beta }) \  x^{\gamma} = x^{\gamma}
$$
The last equality ammounts to
$$
\sum \Sb |\alpha| = d \\  
|\beta| = e \\ \alpha + \beta = \gamma \endSb \frac{d!}{\alpha !} \frac{e!}{\beta !} 
= \frac{(d+e)!}{(\alpha + \beta) !}
$$
To check this formula, let us  multiply 
$$(\sum_{i=1}^n x_i)^d = \sum_{|\alpha| = d} \frac{d!}{\alpha !} \  x^{\alpha} \ \ 
\text{ and } \ \ 
(\sum_{i=1}^n x_i)^e = \sum_{|\beta| = e} \frac{e!}{\beta !} \  x^{\beta}$$
to get
$$\sum_{|\gamma| = d+e} \frac{(d+e)!}{\gamma !} \  x^{\gamma} =
(\sum_{i=1}^n x_i)^{d+e} = 
\sum \Sb |\alpha| = d \\  |\beta| = e \endSb
 \frac{d!}{\alpha !}  \frac{e!}{\beta !} \  x^{\alpha + \beta} =
\sum_{|\gamma| = d+e}  (\sum \Sb |\alpha| = d \\  |\beta| = e \\ \alpha + \beta = \gamma \endSb
  \frac{d!}{\alpha !}  \frac{e!}{\beta !})  \  x^{\gamma}
$$
Equating coefficients of $x^{\gamma}$ we obtain the claim.
\newline
\newline
Now we specialize (2.6) to the case $d=e$. Since the multiplication of the 
symmetric algebra is commutative, by restriction we obtain $GL(V)$-equivariant 
maps, still denoted $\mu$ and $h$
$$S^{2d}(V^*) @>{h}>> S^2(S^d(V^*))  @>{\mu}>> S^{2d}(V^*)$$
given as in (2.5) and satisfying $\mu \circ h = \text{identity}$.
\newline
\newline
This is the desired explicit definition of $h$ as in (2.2). 
Hence, we obtain the sufficient condition:
if the quadratic form $h(F)$ is positive (resp. non-negative) then 
$F \in S^{2d}(V^*)$ is positive (resp. non-negative).

\enddocument